\documentclass[onefignum,onetabnum]{siamart171218}

\usepackage{lipsum}
\usepackage{amsfonts}
\usepackage{graphicx}
\usepackage{epstopdf}
\usepackage{algorithmic}
\ifpdf
  \DeclareGraphicsExtensions{.eps,.pdf,.png,.jpg}
\else
  \DeclareGraphicsExtensions{.eps}
\fi

\newsiamremark{remark}{Remark}
\newsiamremark{hypothesis}{Hypothesis}
\crefname{hypothesis}{Hypothesis}{Hypotheses}
\newsiamthm{claim}{Claim}

\headers{MRT-LB model based four-level FD scheme for 1D diffusion equation}{Y.X. Lin, N. Hong, B.C. Shi and Z. H. Chai}

\title{A multiple-relaxation-time lattice Boltzmann model based four-level finite-difference scheme for one-dimensional diffusion equation\thanks{Submitted to the editors DATE.
\funding{This work was financially supported by the National Natural Science Foundation of China (Grants No. 12072127 and No. 51836003).}}}

\author{   Yuxin Lin  \thanks{School of Mathematics and Statistics, Huazhong University of Science and Technology, Wuhan, 430074, China. }
\and Ning Hong\thanks{School of General Education, Wuchang University of Technology, Wuhan 430223, China.}
  \and Baochang Shi\thanks{School of Mathematics and Statistics, Huazhong University of Science and Technology, Wuhan, 430074, China;
Hubei Key Laboratory of Engineering Modeling and Scientific Computing, Huazhong University of Science and Technology, Wuhan 430074, China.}
\and Zhenhua Chai\footnotemark[2]  \thanks{ Corresponding author. School of Mathematics and Statistics, Huazhong University of Science and Technology, Wuhan, 430074, China; Hubei Key Laboratory of Engineering Modeling and Scientific Computing, Huazhong University of Science and Technology, Wuhan 430074, China. (\email{hustczh@hust.edu.cn})}  }

\usepackage{amsopn}

\ifpdf
\hypersetup{
   pdftitle={A multiple-relaxation-time lattice Boltzmann model based four-level finite-difference scheme for one-dimensional diffusion equation},
  pdfauthor={Y. X. Lin, N. Hong, B. C. Shi and Z.H. Chai}
}
\fi

\begin{document}

\maketitle

\begin{abstract}
In this paper, we first present a multiple-relaxation-time lattice Boltzmann (MRT-LB) model for one-dimensional diffusion equation where the D1Q3 (three discrete velocities in one-dimensional space) lattice structure is considered. Then through the theoretical analysis, we derive an explicit four-level finite-difference scheme from this MRT-LB model. The results show that the four-level finite-difference scheme is unconditionally stable, and through adjusting the weight coefficient $\omega_{0}$ and the relaxation parameters $s_1$ and $s_2$ corresponding to the first and second moments, it can also have a sixth-order accuracy in space. Finally, we also test the four-level finite-difference scheme through some numerical simulations, and find that the numerical results are consistent with our theoretical analysis.
\end{abstract}



\section{Introduction}

In the past decades, the lattice Boltzmann (LB) method, as one of the mesoscopic numerical approaches, not only has gained a great success in the study of complex flows governed by the Navier-Stokes equations \cite{Benzi1992, Chen1998, Aidun2010, Chen2014, Li2016, Liu2016, Wang2019, Sharma2020, Wolf-Gladrow2000, Succi2001, Sukop2006, Guo2013, Kruger2017}, but also can be considered as a numerical solver to some other types of the partial differential equations, for example, the diffusion equations \cite{Dawson1993, Wolf-gladrow1995, Huber2010, Suga2010, Leemput2007}, convection-diffusion equations \cite{van2000, Ginzburg2005,Rasin2005, Ginzburg2009, Shi2009, Chopard2009, Yoshida2010, Chai2013, Huang2014, Chai2016, Aursjo2017, Li2019}, Poisson equation \cite{Hirabayashi2011, Chai2008a, Chai2019}, Kuramoto-Sivashinsky equations \cite{Lai2009, Otomo2017, Chai2018} and some complex equations \cite{Zhong2006, Succi2007, Palpacelli2007, Shi2009b}. Based on the collision term in the LB method, the LB models can be classified into three basic kinds, i.e., the lattice BGK or single-relaxation-time LB (SRT-LB) model \cite{Qian1992}, two-relaxation-time LB (TRT-LB) model \cite{Ginzburg2005} and the multiple-relaxation-time LB (MRT-LB) model \cite{dHumieres1992}. Actually, the SRT-LB and TRT-LB models are two special cases of the MRT-LB model \cite{Chai2020}, and moreover, the MRT-LB model could be more stable and/or more accurate than SRT-LB and TRT-LB models through adjusting some free relaxation parameters \cite{Lallemand2000, Pan2006, Cui2016, Luo2011}. For these reasons, the MRT-LB model is considered in this work.

To build the relation between the LB method and the macroscopic partial differential equations, some asymptotic analysis approaches are usually adopted \cite{Chai2020}, including the Chapman-Enskog expansion \cite{Chapman1970}, Maxwell iteration \cite{Ikenberry1956, Yong2016}, direct Taylor expansion \cite{Holdych2004, Wagner2006, Dubois2008}, recurrence equation \cite{dHumiere2009, Ginzburg2012}, and asymptotic expansion with diffusive scaling \cite{Junk2005}. With the help of these asymptotic analysis methods, one can determine the expression of macroscopic transport efficient, which is related to the relaxation parameter appeared in the LB model. Although the asymptotic analysis methods can be used to illustrate that the LB method is suitable for the macroscopic partial differential equations, and demonstrate the accuracy of the LB method \cite{Dubois2008, Dubois2009a, Dubois2009b}, the relation between the mesoscopic LB method and the macroscopic partial-differential-equation based numerical schemes (hereafter named macroscopic numerical schemes) is still unclear, and furthermore, can we obtain the macroscopic numerical scheme of the LB method? Actually, once the macroscopic numerical scheme of the LB method is obtained, we can not only gain a better understanding on the LB method through the knowledge already available on the macroscopic numerical scheme, but also perform more further research on constructing the mesoscopic LB models for macroscopic partial differential equations.

In the past years, some efforts have been made on this aspect. Junk \cite{Junk2001} and Inamuro \cite{Inamuro2002} found that when the relaxation parameter is equal to unity, the SRT-LB model would reduce to a special macroscopic two-level finite-difference (FD) scheme, and at the diffusive scaling ($\Delta t\propto \Delta x^{2}$, $\Delta t$ and $\Delta x$ are time step and lattice spacing), the macroscopic numerical scheme for the incompressible Naiver-Stokes equations has a second-order convergence rate in space \cite{Junk2001}. Ancona \cite{Ancona1994} demonstrated that for one-dimensional (1D) convection-diffusion equations, the LB method with D1Q2 lattice structure can be written as the classical DuFort-Frankel scheme \cite{Dufort}, which is an explicit three-level second-order FD scheme.
He et al. \cite{He1997JSP} performed a theoretical analysis on the SRT-LB model for several simple flows, and found that when the flows are assumed to be unidirectional and steady-state, the SRT-LB mode is nothing but a macroscopic second-order FD scheme of simplified incompressible Navier-Stokes equations. Then they also obtained the analytical solutions of these simple flows under some commonly used schemes for nonslip boundary conditions, and demonstrated that usually there is a numerical slip velocity at the solid wall which cannot be eliminated effectively in the SRT-LB model. Following the same way, Guo and his collaborators \cite{Guo2008, Chai2010} conducted some analyses on the MRT-LB model for microscale flows, and obtained the similar results as those reported in Ref. \cite{He1997JSP}, while the numerical slip caused by the bounce-back scheme can be eliminated through adjusting the free relaxation parameter corresponding the third-order moment of the distribution function. d'Humi\`{e}re and Ginzburg \cite{dHumiere2009, Ginzburg2012} also carried out a theoretical analysis on the TRT-LB model with recurrence equations, and found that when the relation $\Lambda^{eo}=1/4$ ($\Lambda^{eo}$ is the so-called \emph{magic} parameter) is satisfied, the TRT-LB model would become a macroscopic three-level FD scheme with a second-order accuracy in space. Li et al. \cite{Li2012} conducted an analysis on the SRT-LB model with D1Q2 lattice structure for 1D Burgers equation, and demonstrated that this LB model is equivalent to an explicit three-level second-order FD scheme. Recently, Cui et al. \cite{Cui2016} also carried out a theoretical analysis on the MRT-LB model for the convection-diffusion equation, and found that for the the unidirectional and steady-state problems, the MRT-LB model is just a macroscopic second-order FD scheme of convection-diffusion equation. From the works mentioned above, one can see that under some special conditions, the LB method is equivalent to a special macroscopic second-order FD scheme for a specified partial differential equation, which is also consistent with the results of some asymptotic analyses on the LB method \cite{Chai2020}. However, through choosing the weight coefficients and relaxation parameter properly, Suga \cite{Suga2010} found that the SRT-LB model with D1Q3 lattice structure could be a macroscopic four-level fourth-order FD scheme for 1D diffusion equation. We noted that this work is only limited to the SRT-LB model, it is still unclear whether the more general MRT-LB model can be written as a macroscopic high-order FD scheme, and additionally, can we obtain a more accurate FD scheme from the MRT-LB model where more degrees of freedom in adjusting the relaxation parameters are included? To answer these questions, in this work we first develop a MRT-LB model for 1D diffusion equation where a free weight coefficient $\omega_{0}$ is introduced, and then based on this MRT-LB model, we obtain an equivalent macroscopic four-level FD scheme. Through some theoretical analyses, we show that the four-level FD scheme is unconditionally stable, and can achieve a sixth-order convergence rate in space.

The rest of the paper is organized as follows. In section \ref{section2}, we presented a MRT-LB model for the 1D diffusion equations where the D1Q3 lattice structure is adopted, and then derived an explicit four-level FD scheme from the MRT LB model. Additionally, it is also shown that the macroscopic numerical schemes of the SRT-LB model, TRT-LB model, regularized-LB model and modified-lattice-kinetic model are just some special cases of that of the MRT-LB model. In section \ref{section3}, we investigated the accuracy of the four-level FD scheme, followed by a stability analysis in section \ref{section4}. In section \ref{section5}, we performed some simulations, and found that under some special conditions, the four-level FD scheme has a sixth-order convergence rate in space, which is also consistent with our theoretical analysis. Finally, some conclusions are given in section \ref{section6}.

\section{The multiple-relaxation-time lattice Boltzmann model based four-level finite-difference scheme for one-dimensional diffusion equation}
\label{section2}
In this section, we first developed a MRT-LB model for 1D diffusion equation with a constant diffusion coefficient, and then presented the details on how to obtain an explicit four-level FD scheme from the MRT-LB model.

\subsection{The multiple-relaxation-time lattice Boltzmann model for one-dimensional diffusion equation}

From the mathematical point of view, the 1D diffusion process of mass and heat can be described by the classical diffusion equation,
\begin{equation}\label{eq2-1}
\frac{\partial \phi}{\partial t}=\kappa\frac{\partial^2 \phi}{\partial x^2}+R,
\end{equation}
where $\phi $ is a scalar variable dependent on the space $x$ and time $t$. $\kappa$ is the diffusion coefficient, $R$ is the source term, and in this work, they are assumed to be two constants.
In the framework of LB method, the diffusion equation (\ref{eq2-1}) can be solved efficiently and accurately by different LB models \cite{Dawson1993, Wolf-gladrow1995, Suga2010, Ginzburg2005, Yoshida2010, Chai2016}, here we only consider the more general MRT-LB model for its accuracy and stability in the study of complex problems \cite{Pan2006,Luo2011,Chai2016a}.

The evolution of MRT-LB model for the diffusion equation (\ref{eq2-1}) can be written as \cite{Chai2016,Cui2016}
\begin{equation}\label{eq2-2}
\begin{aligned}
f_{i} (x+\mathbf{c}_{i}\Delta t,t+\Delta t)&=f_{i} (x,t)-(\mathbf{M}^{-1}\mathbf{SM})_{ik} [f_{k} (x,t)-f_{k} ^{eq} (x,t)]\\
&+\Delta t\big[\mathbf{M}^{-1} \big(\mathbf{I}-\frac{\mathbf{S}}{2}\big) \mathbf{M}\big]_{ik} R_k,    (i=-1,0,1),
\end{aligned}
\end{equation}
where $f_{i}(x,t)$ and $f_{i}^{eq}(x,t)$ are the distribution function and equilibrium distribution at position $x$ and time $t$. In the D1Q3 lattice structure, the discrete velocity $\mathbf{c}_{i}$, the transformation matrix $\mathbf{M}$ and the diagonal relaxation matrix $\mathbf{S}$ can be given by
\begin{subequations}
\begin{equation}\label{eq2-3a}
\mathbf{c}_{i}=\left\{\begin{array}{ll}
-c, & \mathrm{\emph{i}=-1},\\
0, & \mathrm{\emph{i}=0},\\
c, & \mathrm{\emph{i}=1},
\end{array}\right.
\end{equation}
\begin{equation}\label{eq2-3b}
\mathbf{M}=\left(\begin{array}{ccc}
1& 1& 1\\ -c & 0 & c\\c^{2}&  -2c^{2}& c^{2}\end{array}\right),
\end{equation}
\begin{equation}\label{eq2-3c}
\mathbf{S}=\left(\begin{array}{ccc}
s_{0}& 0& 0\\ 0& s_{1}& 0\\ 0& 0& s_{2}\end{array}\right),
\end{equation}
\end{subequations}
where $c=\Delta x/\Delta t$ is the lattice speed with $\Delta x$ and $\Delta t$ being lattice spacing and time step, respectively. The diagonal element $s_{i}$ of the relaxation matrix $\mathbf{S}$ is the relaxation parameter corresponding to $i$th moment of the distribution function $f_{i} (x,t)$, and to make the physical transport coefficient (e.g., diffusion coefficient) positive, it should be located in the range $(0,\ 2)$. $R_{i}$ is the discrete source term, and is defined by
\begin{equation}\label{eq2-4}
R_{i}=\omega_{i} R,
\end{equation}
where $\omega_{i}$ is the weight coefficient.

In the LB method, to derive correct macroscopic diffusion equation (\ref{eq2-1}), the equilibrium distribution should be defined as
\begin{equation}\label{eq2-5}
f_{i} ^{eq} (x,t)=\omega_{i} \phi(x,t),
\end{equation}
which satisfies the following conditions \cite{Chai2016},
\begin{subequations}\label{eq2-6}
\begin{equation}\label{eq2-6a}
\sum_{i} f_{i} ^{eq} (x,\ t)=\phi(x,\ t),
\end{equation}
\begin{equation}\label{eq2-6b}
\sum_{i} \mathbf{c}_i f_{i} ^{eq} (x,\ t)= 0,
\end{equation}
\begin{equation}\label{eq2-6c}
\sum_{i} \mathbf{c}_i\mathbf{c}_i f_{i} ^{eq} (x,\ t)= 2\omega_{1}\phi(x,\ t)c^{2},
\end{equation}
\end{subequations}
where the relation $\omega_{-1}=\omega_{1}$ derived from Eq. (\ref{eq2-6b}) is used to obtain Eq. (\ref{eq2-6c}). If the weight coefficient $\omega_{0}$ is considered as a free parameter, then from Eq. (\ref{eq2-6a}) and the condition $\omega_{-1}=\omega_{1}$, one can obtain
\begin{equation}\label{eq2-7}
\omega_{1}=\omega_{-1}=\frac{1-\omega_{0}}{2},
\end{equation}
where $0< \omega_{0} < 1$ which can be used to ensure that all weight coefficients are larger than zero. In addition, the macroscopic variable $\phi(x,\ t)$ can be calculated by
\begin{equation}\label{eq2-8}
\phi(x,t)= \sum_{i} f_{i} (x,\ t)+\frac{\Delta t}{2}R.
\end{equation}
Through the Chapman-Enskog analysis \cite{Chai2016}, one can correctly recover the diffusion equation (\ref{eq2-1}) from the present MRT-LB model with the following relation between the diffusion coefficient and relaxation parameter $s_{1}$,
\begin{equation}\label{eq2-9}
\kappa= 2\omega_{1}\big(\frac{1}{s_{1}}-\frac{1}{2}\big)\frac{\Delta x^{2}}{\Delta t}.
\end{equation}

\subsection{The multiple-relaxation-time lattice Boltzmann model based explicit four-level finite-difference scheme}

In this part, we will show some details on how to derive the macroscopic numerical scheme from the present MRT-LB model. To simplify the following analysis, the notations $f_{i,j}^n=f_{i}(j\Delta x,\ n\Delta t)$ and $\phi_{j}^n=\phi(j\Delta x,\ n\Delta t)$ are introduced. Through substituting the discrete source term $R_{i}$ [Eq. (\ref{eq2-4})] and equilibrium distribution function [Eq. (\ref{eq2-5})] into the evolution equation (\ref{eq2-2}), we have
\begin{subequations}\label{eq2-10}
\begin{equation}\label{eq2-10a}
f_{-1,j}^{n+1}= f_{-1,j+1}^{n}-\frac{s_1}{2}\big(f_{-1,j+1}^{n}-f_{1,j+1}^{n}\big)+\frac{s_2}{2}f_{0,j+1}^{n}-\frac{\omega_{0} s_2}{2} \phi_{j+1}^{n}+\big(\omega_{-1}+\frac{\omega_{0}s_2}{4}\big)\Delta t R ,
\end{equation}
\begin{equation}\label{eq2-10b}
f_{0,j}^{n+1}= f_{0,j}^{n}-s_2 f_{0,j}^{n}+\omega_{0} s_2 \phi_{j}^{n}+\omega_{0}\big(1-\frac{s_2}{2}\big)\Delta t R,
\end{equation}
\begin{equation}\label{eq2-10c}
f_{1,j}^{n+1}= f_{1,j-1}^{n}+\frac{s_1}{2}\big(f_{-1,j-1}^{n}-f_{1, j-1}^{n}\big)+\frac{s_2}{2}f_{0, j-1}^{n}-\frac{\omega_{0} s_2}{2} \phi_{j-1}^{n}+\big(\omega_{1}+\frac{\omega_{0}s_2}{4}\big)\Delta t R,
\end{equation}
\end{subequations}
where Eq. (\ref{eq2-8}) has been used.
To obtain the macroscopic numerical scheme of the MRT-LB model, the distribution functions appeared in Eq. (\ref{eq2-10}) must be replaced by the macroscopic variable $\phi$ at different grid points and time levels. For this purpose, we first conduct a sum of Eq. (\ref{eq2-10}), and derive the following equation,
\begin{equation}\label{eq2-11}
\begin{aligned}
\phi_j^{n+1} &=f_{-1, j+1}^{n}+f_{{0, j}}^{n}+f_{1, j-1}^{n}-\frac{s_1}{2}\big[\big(f_{-1, j+1}^{n}-f_{1, j+1}^{n}\big)-\big(f_{-1, j-1}^{n}-f_{1, j-1}^{n}\big)\big]\\
&+\frac{s_2}{2}\big(f_{0, j-1}^{n}-2f_{0, j}^{n}+f_{0, j+1}^{n}\big)-\frac{\omega_{0} s_2}{2}\big(\phi_{j+1}^{n}-2 \phi_{j}^{n}+\phi_{j-1}^{n}\big)+\frac{3}{2}\Delta t R.
\end{aligned}
\end{equation}
Then from Eq. (\ref{eq2-10b}) we can obtain
\begin{equation}\label{eq2-12}
-f_{0, j}^{n}+s_2 f_{0, j}^{n}-\omega_{0} s_2 \phi_{j}^{n}=-f_{0, j}^{n+1} +\omega_{0}\big(1-\frac{s_2}{4}\big)\Delta t R.
\end{equation}
Based on Eqs. (\ref{eq2-8}), (\ref{eq2-10a}), (\ref{eq2-10c}) and (\ref{eq2-12}), we have
\begin{equation}\label{eq2-13}
f_{1, j+1}^{n}+f_{-1, j-1}^{n}=-f_{0, j}^{n}+\phi_{j}^{n-1}+\frac{1}{2}\Delta t R.
\end{equation}
In addition, from (\ref{eq2-8}) one can derive
\begin{equation}\label{eq2-14}
\begin{aligned}
f_{-1, j+1}^{n}+f_{0, j}^{n}+f_{1, j-1}^{n}&=-f_{-1, j-1}^{n}-f_{1, j+1}^{n}-f_{0, j-1}^{n}-f_{0, j+1}^n+f_{0, j}^{n}\\
&+\phi_{j+1}^{n}+\phi_{j-1}^{n} -\Delta t R\\
&=-\big(f_{0, j-1}^n-2f_{0, j}^n+f_{0, j+1}^n\big)+\phi_{j+1}^{n}+\phi_{j-1}^{n} -\phi_{j}^{n-1}-\frac{3}{2}\Delta t R,
\end{aligned}
\end{equation}
where Eq. (\ref{eq2-13}) has been adopted. Similarly, from Eqs. (\ref{eq2-8}) and (\ref{eq2-13}) we can also obtain
\begin{equation}
\begin{aligned}\label{eq2-15}
&\big(f_{-1, j+1}^{n}-f_{1, j+1}^{n}\big)-\big(f_{-1, j-1}^{n}-f_{1, j-1}^{n}\big)\\
&=-2\big(f_{1, j+1}^{n}+f_{-1, j-1}^{n}\big)-f_{0, j+1}^{n}-f_{0, j-1}^{n}+\phi_{j+1}^{n}+\phi_{j-1}^{n} -\Delta t R\\
&=-\big(f_{0, j-1}^n-2f_{0, j}^n+f_{0, j+1}^n\big)+\phi_{j+1}^{n}+\phi_{j-1}^{n} -2\phi_{j}^{n-1}-2\Delta t R.
\end{aligned}
\end{equation}
Substituting Eqs. (\ref{eq2-14}) and (\ref{eq2-15}) into Eq. (\ref{eq2-11}) yields
\begin{equation}
\begin{aligned}\label{eq2-16}
\phi_{j}^{n+1}
=&\big(\frac{s_1}{2}+\frac{s_2}{2}-1\big)\big(f_{0, j-1}^n-2f_{0, j}^n+f_{0, j+1}^n\big)+\big(1-\frac{s_1}{2}-\frac{\omega_{0} s_2}{2}\big)\phi_{j+1}^{n}\\
&+\omega_{0}s_2\phi_{j}^{n}+\big(1-\frac{s_1}{2}-\frac{\omega_{0} s_2}{2}\big)\phi_{j-1}^{n}+(s_1-1)\phi_{j}^{n-1}+\Delta t s_1 R.
\end{aligned}
\end{equation}
Now we need to give an evaluation of the first term on the right hand side of Eq. (\ref{eq2-16}). Actually, from Eqs. (\ref{eq2-8}) and (\ref{eq2-10}) we have
\begin{equation}
\begin{aligned}\label{eq2-17}
f_{0, j-1}^{n}-2f_{0, j}^{n}+f_{0, j+1}^{n}=&f_{0, j-1}^n+f_{0, j+1}^n-2(\phi_j^n-f_{-1, j}^n-f_{1, j}^n)+ \Delta t R\\
=&2(f_{-1, j}^n+f_{1, j}^n)+f_{0, j-1}^n+f_{0, j+1}^n-2\phi_j^n+\Delta t R\\
=&f_{0, j+1}^{n-1}+2f_{-1, j+1}^{n-1}+f_{0, j-1}^{n-1}+2f_{1, j-1}^{n-1}\\
-&s_1\big[\big(f_{-1, j+1}^{n-1}-f_{1, j+1}^{n-1}\big)-\big(f_{-1, j-1}^{n-1}-f_{1, j-1}^{n-1}\big)\big]-2\phi_j^n+3\Delta tR.
\end{aligned}
\end{equation}
With the help of Eqs. (\ref{eq2-8}) and (\ref{eq2-13}), one can derive the following equation,
\begin{equation}
\begin{aligned}\label{eq2-18}
f_{0,j+1}^{n-1}&+2f_{-1, j+1}^{n-1}+f_{0, j-1}^{n-1}+2f_{1, j-1}^{n-1}\\
=&-f_{0, j+1}^{n-1}-f_{0, j-1}^{n-1}-2\big(f_{-1, j}^{n-1}+f_{1, j}^{n-1}\big)+2\phi_{j+1}^{n-1}+2\phi_{j-1}^{n-1}-2\Delta t R\\
=&-\big(f_{0, j-1}^{n-1}-2f_{0, j}^{n-1}+f_{0, j+1}^{n-1}\big)+2\phi_{j+1}^{n-1}+2\phi_{j-1}^{n-1}-2\phi_{j}^{n-2}-3\Delta t R.
\end{aligned}
\end{equation}
Substituting Eqs. (\ref{eq2-15}) and (\ref{eq2-18}) into Eq. (\ref{eq2-17}) gives rise to
\begin{equation}
\begin{aligned}\label{eq2-19}
f_{0, j-1}^{n}-2f_{0, j}^{n}+f_{0, j+1}^{n}=&(s_1-1)\big(f_{0,j-1}^{n-1}-2f_{0, j}^{n-1}+f_{0, j+1}^{n-1}\big)\\
-&2\phi_{j}^{n}+(2-s_1)\phi_{j-1}^{n-1}+(2-s_1)\phi_{j}^{n-2}-(2s_1-2)\phi_{j}^{n-2}\\
+&2s_1\Delta t R.
\end{aligned}
\end{equation}
If we substitute Eq. (\ref{eq2-19}) into Eq. (\ref{eq2-16}), one can obtain
\begin{equation}
\begin{aligned}\label{eq2-20}
\phi_{j}^{n+1}=&(s_1-1)\big(\frac{s_1}{2}+\frac{s_2}{2}-1\big)\big(f_{{0}_{j-1}}^{n-1}-2f_{{0}_{j}}^{n-1}+f_{{0}_{j+1}}^{n-1}\big)\\
+&\big(1-\frac{s_1}{2}-\frac{\omega_{0}s_2}{2}\big)\phi_{j-1}^n+[\omega_{0} s_2-2\big(\frac{s_1}{2}+\frac{s_2}{2}-1\big)\big]\phi_{j}^n\\
+&\big(1-\frac{s_1}{2}-\frac{\omega_{0}s_2}{2}\big)\phi_{j+1}^n+(2-s_1)\big(\frac{s_1}{2}+\frac{s_2}{2}-1\big)\phi_{j-1}^{n-1}+(s_1-1)\phi_{j}^{n-1}\\
+&(2-s_1)\big(\frac{s_1}{2}+\frac{s_2}{2}-1\big)\phi_{j+1}^{n-1}+(2s_1-2)\big(\frac{s_1}{2}+\frac{s_2}{2}-1\big)\phi_{j}^{n-2}\\
+&s_1(s_1+s_2-1)\Delta t R.
\end{aligned}
\end{equation}
With the aid of Eq. (\ref{eq2-16}), we can rewrite Eq. (\ref{eq2-20}) as
\begin{equation}
\begin{aligned}\label{eq2-21a}
\phi_{j}^{n+1}=&\big(1-\frac{s_1}{2}-\frac{\omega_{0} s_2}{2}\big)\phi_{j-1}^n+[(\omega_{0}-1)s_2+1]\phi_j^n+\big(1-\frac{s_1}{2}-\frac{\omega_{0}s_2}{2}\big)\phi_{j+1}^n\\
+&\big(\frac{\omega_{0} s_1 s_2}{2}-\frac{s_1 s_2}{2}-\frac{\omega_{0} s_2}{2}+\frac{s_1}{2}+s_2-1\big)\phi_{j-1}^{n-1}\\
+&(-\omega_{0} s_1 s_2+\omega_{0} s_2+s_1-1)\phi_{j}^{n-1}\\
+&\big(\frac{\omega_{0}s_1s_2}{2}-\frac{s_1s_2}{2}-\frac{\omega_{0}s_2}{2}+\frac{s_1}{2}+s_2-1\big)\phi_{j+1}^{n-1}\\
+&(s_1-1)(s_2-1)\phi_{j}^{n-2}+s_1 s_2\Delta t R,
\end{aligned}
\end{equation}
which can also be written as
\begin{equation}
\begin{aligned}\label{eq2-21}
\phi_{j}^{n+1}=&\alpha_{1}\phi_{j-1}^n+\alpha_{2}\phi_j^n+\alpha_{1}\phi_{j+1}^n\\
+&\beta_{1}\phi_{j-1}^{n-1}+\beta_{2}\phi_{j}^{n-1}+\beta_{1}\phi_{j+1}^{n-1}\\
+&\gamma\phi_{j}^{n-2}+\delta\Delta t R,
\end{aligned}
\end{equation}
where the parameters $\alpha_{i}$ ($i=1, 2$), $\beta_{i}$, $\gamma$ and $\delta$ are given by
\begin{equation}\label{eq2-21c}
\begin{aligned}
\alpha_1&=1-\frac{s_1}{2}-\frac{\omega_{0} s_2}{2}, & \alpha_2&=(\omega_{0}-1)s_2+1, \\
\beta_1&=\frac{\omega_{0} s_1 s_2}{2}-\frac{s_1 s_2}{2}-\frac{\omega_{0} s_2}{2}+\frac{s_1}{2}+s_2-1, & \beta_2&=-\omega_{0} s_1 s_2+\omega_{0} s_2+s_1-1,\\
\gamma &=(s_1-1)(s_2-1), & \delta&=s_1 s_2.
\end{aligned}
\end{equation}
Here we would like to point out that Eq. (\ref{eq2-21}) is the the mesoscopic MRT-LB model based macroscopic explicit four-level FD scheme for 1D diffusion equation. \\
\textbf{Remark I:} If $s_0=s_1=s_2=\omega$ ($\omega$ is the relaxation parameter in the SRT-LB model) and $R=0$, we can obtain the following SRT-LB model based four-level FD scheme from Eq. (\ref{eq2-21}),
\begin{equation}
\begin{aligned}\label{eq2-22}
\phi_{j}^{n+1}=&\alpha_1\phi_{j-1}^n+\alpha_2\phi_{j}^{n}+\alpha_1\phi_{j+1}^{n}\\	
+&\beta_1\phi_{j-1}^{n-1}+\beta_2\phi_{j}^{n-1}+\beta_1\phi_{j+1}^{n-1}\\
+&\gamma\phi_{j}^{n-2},
\end{aligned}
\end{equation}
with the following parameters,
\begin{equation}
\begin{aligned}
\alpha_1=\Omega+\omega_{1}\omega,\ \ & \alpha_2=\Omega+(1-2\omega_{1})\omega, \\
\beta_1=-[\Omega+(1-\omega_{1})\omega]\Omega,\ \ \beta_2&=-(\Omega+2\omega_{1}\omega)\Omega,\ \ \gamma&=\Omega^{2},
\end{aligned}
\end{equation}
where $\Omega=1-\omega$. It is clear that the SRT-LB model based macroscopic numerical scheme (\ref{eq2-22}), as a special case of the present MRT-LB model based four-level scheme, is the same as that reported in the previous work \cite{Suga2010}.\\
\textbf{Remark II:} If $s_0=s_2=s^{+}$ and $s_{1}=s^{-}$ with $s^{+}$ and $s^{-}$ denoting the relaxation parameters corresponding to the symmetric and antisymmetric modes \cite{Ginzburg2005}, we can obtain the TRT-LB model based four-level FD scheme from Eq. (\ref{eq2-21}). However, when $s_{0}=s_{2}=1$ and $s_{1}=\omega$, one can derive the regularized-LB model \cite{Latt2006,Wang2015} based macroscopic numerical scheme from Eq. (\ref{eq2-21}). In addition, if $s_0=s_2=\omega$ and $s_{1}=\omega/(1-\omega\eta)$ with $\eta$ being an adjusting parameter, we can also obtain the modified-lattice-kinetic model \cite{Yang2014, Wang2015a, Wang2018} based macroscopic numerical scheme from Eq. (\ref{eq2-21}).

\section{The accuracy analysis of the multiple-relaxation-time lattice Boltzmann model based macroscopic numerical scheme}\label{section3}
We now performed an accuracy analysis on the macroscopic four-level FD scheme (\ref{eq2-21}). To do this, we first conducted the Taylor expansion to Eq. (\ref{eq2-21}) at the position $x=j\Delta x$ and time $t=n\Delta t$, and after some algebraic manipulations, one can obtained
\begin{equation}
\begin{aligned}\label{eq3-1}
(1+2\beta_1+\beta_2+2\gamma)\big[\frac{\partial \phi}{\partial t}\big]_j^n=&(\alpha_1+\beta_1)\frac{\Delta x^2}{\Delta t}\big[\frac{\partial^2 \phi}{\partial x^2}\big]_j^n+\delta R\\
+&\frac{1}{12}(\alpha_1+\beta_1)\frac{\Delta x^4}{\Delta t}\big[\frac{\partial^4 \phi}{\partial x^4}\big]_j^n+\frac{1}{360}(\alpha_1+\beta_1)\frac{\Delta x^6}{\Delta t}\big[\frac{\partial^6 \phi}{\partial x^6}\big]_j^n\\
-&\beta_1\Delta x^2\big[\frac{\partial^3 \phi}{\partial x^2\partial t}\big]_j^n-\frac{\beta_1}{12}\Delta x^4\big[\frac{\partial^5 \phi}{\partial x^4\partial t}\big]_j^n\\
+&\frac{\beta_1}{2}\Delta x^2\Delta t\big[\frac{\partial^4 \phi}{\partial x^2\partial t^2}\big]_j^n+\frac{1}{2}\big(2\beta_1+\beta_2+4\gamma-1\big)\Delta t\big[\frac{\partial^2 \phi}{ \partial t^2}\big]_j^n\\
-&\frac{1}{6}\big(2\beta_1+\beta_2+8\gamma+1\big)\Delta t^2\big[\frac{\partial^3 \phi}{ \partial t^3}\big]_j^n+\cdots.
\end{aligned}
\end{equation}
Substituting Eq. (\ref{eq2-21c}) into above equation yields
\begin{equation}
\begin{aligned}\label{eq3-2}
\big[\frac{\partial \phi}{\partial t}\big]_j^n=&\kappa\big[\frac{\partial^2 \phi}{\partial x^2}\big]_j^n+R\\
+&\frac{1}{12}\kappa\Delta x^2\big[\frac{\partial^4 \phi}{\partial x^4}\big]_j^n+\frac{1}{360}\kappa\Delta x^4\big[\frac{\partial^6 \phi}{\partial x^6}\big]_j^n\\
-&\big(\frac{\omega_{0}}{2}-\frac{1}{2}-\frac{\omega_{0}}{2s_1}+\frac{1}{2s_2}+\frac{1}{s_1}-\frac{1}{s_1s_2}\big)\Delta x^2\big[\frac{\partial^3 \phi}{\partial x^2\partial t}\big]_j^n\\
-&\big(\frac{\omega_{0}}{24}-\frac{1}{24}-\frac{\omega_{0}}{24s_{1}}+\frac{1}{24s_2}+\frac{1}{12s_1}-\frac{1}{12s_1s_2}\big)\Delta x^4\big[\frac{\partial^5 \phi}{\partial x^4\partial t}\big]_j^n\\
+&\big(\frac{\omega_{0}}{4}-\frac{1}{4}-\frac{\omega_{0}}{4s_1}+\frac{1}{4s_2}+\frac{1}{2s_1}-\frac{1}{2s_1s_2}\big)\Delta x^2\Delta t\big[\frac{\partial^4 \phi}{\partial x^2\partial t^2}\big]_j^n\\
+&\big(\frac{3}{2}-\frac{1}{s_1}-\frac{1}{s_2}\big)\Delta t\big[\frac{\partial^2 \phi}{ \partial t^2}\big]_j^n+\big(\frac{1}{s_2}+\frac{1}{s_1}-\frac{7}{6}-\frac{1}{s_1s_2}\big)\Delta t^2\big[\frac{\partial^3 \phi}{ \partial t^3}\big]_j^n+\cdots,
\end{aligned}
\end{equation}
where Eq. (\ref{eq2-9}) has been used to derive above equation.

According to the diffusion equation (\ref{eq2-1}), we have the following relations,
\begin{equation}
\begin{aligned}\label{eq3-3}
\big[\frac{\partial^2 \phi}{ \partial t^2}\big]_j^n =\kappa^2\big[\frac{\partial^4 \phi}{\partial x^4}\big]_j^n,\ \ & \big[\frac{\partial^3 \phi}{ \partial t^3}\big]_j^n=\kappa^3\big[\frac{\partial^6 \phi}{\partial x^6}\big]_j^n,\\
\big[\frac{\partial^5 \phi}{\partial x^4\partial t}\big]_j^n =\kappa\big[\frac{\partial^6 \phi}{\partial x^6}\big]_j^n,\ \  \big[\frac{\partial^3 \phi}{\partial x^2\partial t}\big]_j^n & =\kappa\big[\frac{\partial^4 \phi}{\partial x^4}\big]_j^n,\ \ \big[\frac{\partial^4 \phi}{\partial x^2\partial t^2}\big]_j^n =\kappa^2\big[\frac{\partial^6 \phi}{\partial x^6}\big]_j^n,
\end{aligned}
\end{equation}
Substituting Eq. (\ref{eq3-3}) into Eq. (\ref{eq3-2}), we have
\begin{equation}
\begin{aligned}\label{eq3-4}
\big[\frac{\partial \phi}{\partial t}\big]_j^n=&\kappa\big[\frac{\partial^2 \phi}{\partial x^2}\big]_j^n+R\\
+&\frac{1}{s_1s_2}\big[\frac{s_1s_2}{12}-\big(\frac{\omega_{0}s_1s_2}{2}-\frac{s_1s_2}{2}-\frac{\omega_{0}s_2}{2}+\frac{s_1}{2}+s_2-1\big)\\
+&\big(\frac{3s_1s_2}{2}-s_2-s_1\big)\epsilon\big]\kappa\Delta x^2\big[\frac{\partial^4 \phi}{\partial x^4}\big]_j^n\\
+&\frac{1}{s_1s_2}\big[\frac{s_1s_2}{360}-\big(\frac{\omega_{0}s_1s_2}{24}-\frac{s_1s_2}{24}-\frac{\omega_{0}s_2}{24}+\frac{s_1}{24}+\frac{s_2}{12}-\frac{1}{12}\big)\\
+&\big(\frac{\omega_{0}s_1s_2}{4}-\frac{s_1s_2}{4}-\frac{\omega_{0}s_2}{4}+\frac{s_1}{4}+\frac{s_2}{2}-\frac{1}{2}\big)\epsilon\\
+&\big(-\frac{7s_1s_2}{6}+s_2+s_1-1\big)\epsilon^2\big]\kappa\Delta x^4\big[\frac{\partial^6 \phi}{\partial x^6}\big]_j^n\\
+&O(\Delta x^6+\Delta t^3),
\end{aligned}
\end{equation}
where $\epsilon$ is the mesh Fourier number, and is defined by
\begin{equation}\label{eq3-5}
\epsilon=\frac{\kappa\Delta t}{\Delta x^{2}}=2\omega_{1}\big(\frac{1}{s_{1}}-\frac{1}{2}\big)=(1-\omega_{0})\big(\frac{1}{s_{1}}-\frac{1}{2}\big).
\end{equation}
From above discussion, it is clear that for a given diffusion coefficient $\kappa$ or mesh Fourier number $\epsilon$, one can obtain an explicit four-level FD scheme with the third-order accuracy in time and sixth-order accuracy in space once the following conditions are satisfied, i.e., the second and fourth-order truncation errors in Eq. (\ref{eq3-4}) are equal to zero,
\begin{subequations}\label{eq3-6}
\begin{equation}\label{eq3-6a}
\frac{s_1s_2}{12}-\big(\frac{\omega_{0}s_2}{2}+\frac{s_1}{2}-1\big)+\big(\frac{s_1s_2}{2}-s_2-s_1\big)\epsilon=0\\
\end{equation}
\begin{equation}\label{eq3-6b}
\frac{s_1s_2}{360}-\frac{1}{12}\big(\frac{\omega_{0}s_2}{2}+\frac{s_1}{2}-1\big)-\frac{1}{2}\big(\frac{s_1s_2}{6}-\frac{\omega_{0}s_2}{2}-\frac{s_1}{2}+1\big)\epsilon
+\big(-\frac{2s_1s_2}{3}+s_2+s_1-1\big)\epsilon^2=0,
\end{equation}
\end{subequations}
where Eq. (\ref{eq3-5}) has been applied. We note that compared to the SRT-LB model based four-level FD scheme for 1D equation \cite{Suga2010}, the present MRT-LB model based macroscopic numerical scheme can be more accurate through adjusting the weight coefficient $\omega_{0}$ and relaxation parameters $s_{1}$ and $s_{2}$ to satisfy Eq. (\ref{eq3-6}). In the following, some remarks on the MRT-LB model based macroscopic numerical scheme are listed.\\
\textbf{Remark I:} It is found that the relaxation parameter $s_{0}$ corresponding to the zeroth moment of distribution function $f_{i}(x, t)$ (or conservative variable $\phi$) does not appear in the macroscopic numerical scheme (\ref{eq2-21}), and thus it has no influence on the numerical results. This also explains why the relaxation parameter $s_{0}$ in the MRT-LB method can be chosen arbitrarily. However, unlike the relaxation parameter $s_{0}$, the relaxation parameter $s_{2}$ corresponding to the second-order moment of distribution function has an important influence on the macroscopic numerical scheme (\ref{eq2-21}) [see Eq. (\ref{eq3-6})], and also the numerical results. Actually, these results have also been reported in some previous works on the MRT-LB model for (convection) diffusion equations \cite{Chai2016, Cui2016, Chai2016a}.\\
\textbf{Remark II:} When the conditions of Eqs. (\ref{eq3-5}) and (\ref{eq3-6}) are satisfied, the MRT-LB model based macroscopic four-level numerical scheme (\ref{eq2-21}) has a sixth-order accuracy in space. However, if only Eqs. (\ref{eq3-5}) and (\ref{eq3-6a}) hold, the macroscopic four-level numerical scheme would have a fourth-order accuracy in space. Compared to the SRT-LB model based macroscopic four-level fourth-order numerical scheme \cite{Suga2010}, the present macroscopic numerical scheme has another distinct characteristic, i.e., there is a free relaxation parameter $s_{2}$, which can be used to remove the discrete effect of anti-bounce-back scheme for the Dirichlet boundary conditions \cite{Cui2016, Chai2016a}.\\
\textbf{Remark III:} Compared to the classical two-level LB method, the implementation of the macroscopic four-level numerical scheme needs the initial values of variable $\phi$ at first three time levels, and usually to complete the initialization, we must adopt some other numerical schemes to obtain the values of variable $\phi$ at the second and third time levels. In addition, the inclusion of the extra time levels also brings a larger memory requirement to store the variable $\phi$.

\begin{figure}[htb]
\centering
\includegraphics[width=0.8\textwidth]{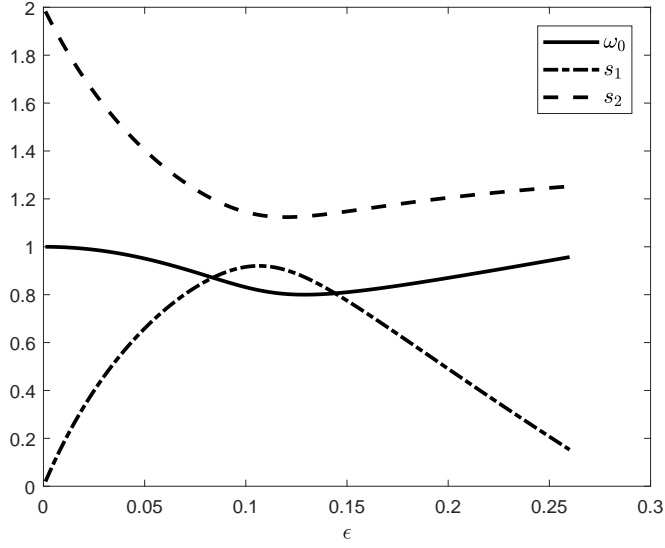}
\caption{The weight coefficient $\omega_{0}$ and the relaxation parameters $s_1$ and $s_2$ as a function of the discretization parameter $\epsilon$.}
\label{fig1}
\end{figure}

Finally, to implement the MRT-LB model based macroscopic sixth-order numerical scheme (\ref{eq2-21}), we must determine the weight coefficient $\omega_{0}$ and relaxation parameters $s_{1}$ and $s_{2}$ from Eqs. (\ref{eq3-5}) and (\ref{eq3-6}) for a specified $\epsilon$, while due to the coupling and nonlinearity of these equations, it is difficult to derive the explicit expressions of $\omega_{0}$, $s_{1}$ and $s_{2}$ in terms of $\epsilon$. For this reason, we numerically solved Eqs. (\ref{eq3-5}) and (\ref{eq3-6}), and plotted $\omega_{0}$, $s_1$ and $s_2$ as a function of the mesh Fourier number $\epsilon$ in Fig. \ref{fig1}. Here it should be noted that only under the condition of $0< \epsilon\le \epsilon_{max}$ with $\epsilon_{max}$ being about 0.262, we can obtain the real roots of Eqs. (\ref{eq3-5}) and (\ref{eq3-6}). As seen from Fig. \ref{fig1}, for a given $\epsilon$, the weight coefficient $\omega_{0}$ and relaxation parameters $s_{1}$ and $s_{2}$ are located in the ranges of $0.8\le \omega_{0}<1$, $0 < s_1 \le s_{1,max}$ and $s_{2,min}\le s_2 < 2.0$ with $s_{1,max}$ and $s_{2,min}$ being very close to 0.920 and 1.124.

\section{The stability analysis of the multiple-relaxation-time lattice Boltzmann model based macroscopic numerical scheme}\label{section4}

In this section, we will prove that the MRT-LB model based four-level FD scheme (\ref{eq2-21}) is unconditionally stable. To this end, we first neglected the source term $R$ and replaced the $\phi_{j}^{n}$ in Eq. (\ref{eq2-10}) by the distribution function $f_{i, j}^{n}$ through the relation (\ref{eq2-8}), and then take the discrete Fourier transform of $f_{i, j}^n$ in Eq. (\ref{eq2-10}) to obtain the following matrix equation \cite{Thomas1995},
\begin{equation}\label{eq4-1}
\widehat{\mathbf{U}}_j^{n+1}=\mathbf{G}(\theta, \omega_{0}, s_1,s_2)\widehat{\mathbf{U}}_j^{n},
\end{equation}
where $\widehat{\mathbf{U}}_j^{n}$ is the discrete Fourier transform of $f_{i, j}^n$ $(i=-1,0,1)$. $\mathbf{G}$ is the amplification matrix of the scheme, and is given by
\begin{equation}\label{eq4-2}
\mathbf{G}=\left(\begin{array}{cccc}
\big(1-\frac{s_1}{2}-\frac{\omega_{0} s_2}{2}\big)e^{i\theta}& \big(\frac{s_2}{2}-\frac{\omega_{0} s_2}{2}\big)e^{i\theta}& \big(\frac{s_1}{2}-\frac{\omega_{0} s_2}{2}\big)e^{i\theta}\\
\omega_{0} s_2& \omega_{0}s_2-s_2+1& \omega_{0} s_2\\
\big(\frac{s_1}{2}-\frac{\omega_{0}s_2}{2}\big)e^{-i\theta}& \big(\frac{s_2}{2}-\frac{\omega_{0} s_2}{2}\big)e^{-i\theta}& \big(1-\frac{s_1}{2}-\frac{\omega_{0} s_2}{2}\big)e^{-i\theta}\end{array}\right),
\end{equation}
where $-\pi \leq \theta \leq \pi$. On the other hand, if we take the discrete Fourier transform of Eq. (\ref{eq2-21}), one can obtain the amplification matrix $\mathbf{H}$ of the macroscopic numerical scheme,
\begin{equation}\label{eq4-3}
\mathbf{H}=\left(\begin{array}{cccc}
2\alpha_{1}\cos\theta + \alpha_{2} & 2\beta_{1}\cos\theta + \beta_{2} & \gamma \\
1 & 0 & 0 \\
0 & 1 & 0 \end{array}\right).
\end{equation}
Although the amplification matrix $\mathbf{H}$ is different from $\mathbf{G}$, due to the equivalence between the MRT-LB model [Eq. (\ref{eq2-2}) or (\ref{eq2-10})] and the macroscopic numerical scheme (\ref{eq2-21}), the characteristic polynomials of them are identical, and can be expressed as
\begin{equation}\label{eq4-4}
p(\lambda)=\lambda^{3}+p_{2}\lambda^{2}+p_{1}\lambda+p_{0},
\end{equation}
where the coefficients $p_{0}$, $p_{1}$ and $p_{2}$ are given by
\begin{equation}\label{eq4-4a}
\begin{aligned}
p_0 & =(s_{1}-1)(1-s_{2}),\\
p_1 & =(s_{1}-1)(s_{2}\omega_{0}-1)+[(s_{1}-2)(s_{2}-1)+s_{2}\omega_{0}(1-s_{1})]\cos\theta,\\
p_2 & = s_{2}-s_{2}\omega_{0}-1+(s_{2}\omega_{0}+s_{1}-2)\cos\theta.
\end{aligned}
\end{equation}
In the following, we would show that the roots of the characteristic polynomial $p(\lambda)$ denoted by $\lambda_{k}$ ($k$=1, 2 and 3) satisfy the condition $|\lambda_{k}|\leq 1$.

With the following linear fractional transformation,
\begin{equation}\label{eq4-5}
\lambda=\frac{1+z}{1-z}, \ z\in \mathcal{C},
\end{equation}
the unit circle $|\lambda|=1$ and the field $|\lambda|<1$ are mapped to the imaginary axis [Re$(z)=0$] and left-half plane [Re$(z)<0$], and vise verse. Here $\mathcal{C}$ and Re denote the complex-number field and the real part of a complex number. Substituting Eq. (\ref{eq4-5}) into Eq. (\ref{eq4-4}), we have
\begin{equation}\label{eq4-6}
\begin{aligned}
(1-z)^{3}p\big(\frac{1+z}{1-z}\big)& =(1+z)^{3}+p_{2}(1-z)(1+z)^{2}+p_{1}(1-z)^{2}(1+z)+p_{0}(1-z)^{3}\\
& = (1-p_{0}+p_{1}-p_{2})z^{3}+(3+3p_0-p_1-p_2)z^{2}\\
& +(3-3p_0-p_1+p_2)z+(1+p_{0}+p_{1}+p_{2}).
\end{aligned}
\end{equation}
To ensure that the roots of characteristic polynomial $p(\lambda)$ are located in the field $|\lambda|<1$, the following conditions must be satisfied, i.e., the Routh-Hurwitz stability criterion \cite{Routh1887, Hurwitz1895, Gantmacher1959},
\begin{subequations}\label{eq4-7}
\begin{equation}\label{eq4-7a}
1-p_{0}+p_{1}-p_{2}>0
\end{equation}
\begin{equation}\label{eq4-7b}
3+3p_0-p_1-p_2>0,
\end{equation}
\begin{equation}\label{eq4-7c}
3-3p_0-p_1+p_2> 0,
\end{equation}
\begin{equation}\label{eq4-7d}
1+p_{0}+p_{1}+p_{2}>0,
\end{equation}
\begin{equation}\label{eq4-7e}
1-p_1+p_0p_2-p_0^{2}>0.
\end{equation}
\end{subequations}
After the sums of Eqs. (\ref{eq4-7a}) and (\ref{eq4-7c}), and Eqs. (\ref{eq4-7b}) and (\ref{eq4-7d}), we can equivalently rewrite Eq. (\ref{eq4-7}) as
\begin{subequations}\label{eq4-8}
\begin{equation}\label{eq4-8a}
1-p_{0}+p_{1}-p_{2}>0
\end{equation}
\begin{equation}\label{eq4-8b}
1-p_0>0,
\end{equation}
\begin{equation}\label{eq4-8c}
1+p_0> 0,
\end{equation}
\begin{equation}\label{eq4-8d}
1+p_{0}+p_{1}+p_{2}>0,
\end{equation}
\begin{equation}\label{eq4-8e}
1-p_1+p_0p_2-p_0^{2}>0.
\end{equation}
\end{subequations}

Actually, under the condition of $\cos\theta\neq 1$, we can first obtain
\begin{subequations}\label{eq4-9}
\begin{equation}\label{eq4-9a}
\begin{aligned}
1-p_{0}+p_{1}-p_{2} = (2-s_1)(2-s_2)(1+\cos\theta)+s_1s_2\omega_{0}(1-\cos\theta)>0,
\end{aligned}
\end{equation}
\begin{equation}\label{eq4-9b}
1+p_{0}=1+(s_{1}-1)(1-s_{2})>0,
\end{equation}
\begin{equation}\label{eq4-9c}
1-p_{0}=1-(s_{1}-1)(1-s_{2})>0,
\end{equation}
\begin{equation}\label{eq4-9d}
1+p_{0}+p_{1}+p_{2}=s_{2}(1-\cos\theta)(2-s_{1})(1-\omega_{0})> 0,
\end{equation}
\end{subequations}
where $0< \omega_{0} <1$ and $0< s_{1}, s_{2} <2$ have been used.

To prove Eq. (\ref{eq4-8e}), we first introduce the parameters $A$ and $B$,
\begin{equation}\label{eq4-10}
\begin{aligned}
A & =s_1(1-s_2)(2-s_1)+\omega_{0}s_2(1-s_1)(2-s_2), \\
B & =s_1s_2(s_1+s_2-s_1s_2)=s_1s_2[1-(1-s_1)(1-s_2)]>0,
\end{aligned}
\end{equation}
and can express the part on the left hand side of Eq. (\ref{eq4-8e}) as
\begin{equation}\label{eq4-9e}
\begin{aligned}
1-p_{0}^{2}-p_{1}+p_{0}p_{2} & =[s_1(1-s_2)(2-s_1)+\omega_{0}s_2(1-s_1)(2-s_2)](1-\cos\theta)\\
& +s_1s_2(s_1+s_2-s_1s_2)=A(1-\cos\theta)+B.
\end{aligned}
\end{equation}
(i) If $0<s_1\leq 1$ and $0<s_2\leq 1$, we can obtain $A\geq 0$ and the following equation,
\begin{equation}\label{eq4-10e}
\begin{aligned}
1-p_{0}^{2}-p_{1}+p_{0}p_{2} =A(1-\cos\theta)+B>0.
\end{aligned}
\end{equation}
(ii) If $0<s_1\leq 1$ and $1<s_2<2$, we have
\begin{equation}\label{eq4-11}
A=s_1(1-s_2)(2-s_1)+\omega_{0}s_2(1-s_1)(2-s_2)\geq s_1(1-s_2)(2-s_1),
\end{equation}
which can be used to derive
\begin{equation}\label{eq4-11a}
\begin{aligned}
1-p_{0}^{2}-p_{1}+p_{0}p_{2}\geq 2s_1(1-s_2)(2-s_1)+B=s_1(2-s_2)[1+(1-s_1)(1-s_2)]>0.
\end{aligned}
\end{equation}
(iii) If $1<s_1< 2$ and $0<s_2\leq 1$, one can obtain
\begin{equation}\label{eq4-12}
A=s_1(1-s_2)(2-s_1)+\omega_{0}s_2(1-s_1)(2-s_2)> s_2(1-s_1)(2-s_2),
\end{equation}
then we have
\begin{equation}\label{eq4-12a}
\begin{aligned}
1-p_{0}^{2}-p_{1}+p_{0}p_{2}> 2s_2(1-s_1)(2-s_2)+B=s_2(2-s_1)[1+(1-s_1)(1-s_2)]>0.
\end{aligned}
\end{equation}
(iv) If $1<s_1< 2$ and $1<s_2<2$, one can derive
\begin{equation}\label{eq4-13}
A=s_1(1-s_2)(2-s_1)+\omega_{0}s_2(1-s_1)(2-s_2)>s_1(1-s_2)(2-s_1)+s_2(1-s_1)(2-s_2).
\end{equation}
With the help of above equation and let $C=s_1(1-s_2)(2-s_1)+s_2(1-s_1)(2-s_2)$, we can obtian
\begin{equation}\label{eq4-13a}
\begin{aligned}
1-p_{0}^{2}-p_{1}+p_{0}p_{2} & > 2A+B> 2C+B\\
&=(2-s_1)(2-s_2)[1+(s_1-1)(1-s_2)]>0.
\end{aligned}
\end{equation}
Based on above results (i)-(iv), one can find that Eq. (\ref{eq4-8e}) indeed holds under the conditions of $0<\omega_{0}<1$ and $0<s_1, s_2<2$. Thus the roots of characteristic polynomial $p(\lambda)$ are located in the field $|\lambda|<1$ under the condition of $\cos\theta\neq 1$.

Now let us focus on the case of $\cos\theta=1$. We note that the roots of the characteristic polynomial $p(\lambda)$ are continuous functions of $\cos\theta$, and hence the roots of characteristic polynomial satisfy the condition $|\lambda_{k}|\leq1$ ($k$=1, 2 and 3). In addition, we would also like to point out that for the special case of $\cos\theta=1$, one can adopt the reductive approach \cite{Miller} to obtain $|\lambda_{k}|\leq1$ ($k$=1, 2 and 3).

From above discussion, we can find that the roots of characteristic polynomial satisfy the condition $|\lambda_{k}|\leq1$ ($k$=1, 2 and 3), thus the present MRT-LB model based macroscopic numerical scheme (\ref{eq2-21}) is unconditionally stable.

\section{Numerical results and discussion}\label{section5}
\begin{figure}
\centering
\includegraphics[width=0.8\textwidth]{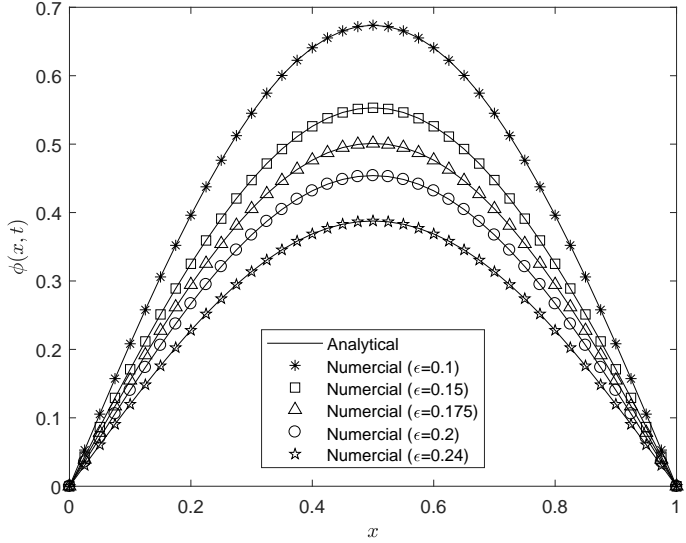}
\caption{The numerical and analytical solutions under different values of discretization parameter $\epsilon$.}
\label{fig2}
\end{figure}

To test the capacity of the MRT-LB model based macroscopic four-level numerical scheme (\ref{eq2-21}), we considered the diffusion equation (\ref{eq2-1}) with the following initial and boundary conditions \cite{Fletcher1988},
\begin{equation}\label{eq5-1}
\begin{aligned}
\phi(x,0)=\sin(\pi x), ~0\le x\le 1,\\
\phi(0,t)=\phi(1,t)=0, ~t>0,
\end{aligned}
\end{equation}
and obtained the analytical solution of this problem,
\begin{equation}\label{eq5-2}
\phi(x,t)=\sin(\pi x) e^{-\kappa \pi^2 t}.
\end{equation}
In the implementation of the macroscopic four-level numerical scheme (\ref{eq2-21}), the analytical solution (\ref{eq5-2}) is used to initialize the variable $\phi$ at first three time levels.
We first performed some simulations under different values of the mesh Fourier number $\epsilon$ (or different diffusion coefficients for the specified time and space steps), and plotted the results in Fig. \ref{fig2} where $\Delta x=0.025$, $\Delta t=0.01875$, the weight coefficient $\omega_{0}$ and the relaxation parameters $s_{1}$ and $s_{2}$ are determined from Eqs. (\ref{eq3-5}) and (\ref{eq3-6}). As shown in this figure, the numerical results are in good agreement with the analytical solutions at time $t=12$.

In addition, to quantitatively measure the deviation between the numerical results and analytical solutions, the following root-mean-square error ($RMSE$) is adopted \cite{Suga2010},
\begin{equation}\label{eq5-3}
RMSE=\sqrt{\frac{\sum_{j=1}^{N}\big[\phi(j\Delta x, n\Delta t)- \phi^{*} (j\Delta x, n\Delta t)\big]^2}{N}},
\end{equation}
where $N$ is the number of grid points, $\phi$ and $\phi^{*}$ are the numerical and analytical solutions. Based on the definition of $RMSE$, one can evaluate the convergence rate ($C_R$) of numerical scheme with the following formula,
\begin{equation}\label{40}
C_R=\log\big(RMSE_{\Delta x}/RMSE_{\Delta x/2}\big)/\log2.
\end{equation}

\begin{figure}
\centering
\includegraphics[width=0.8\textwidth]{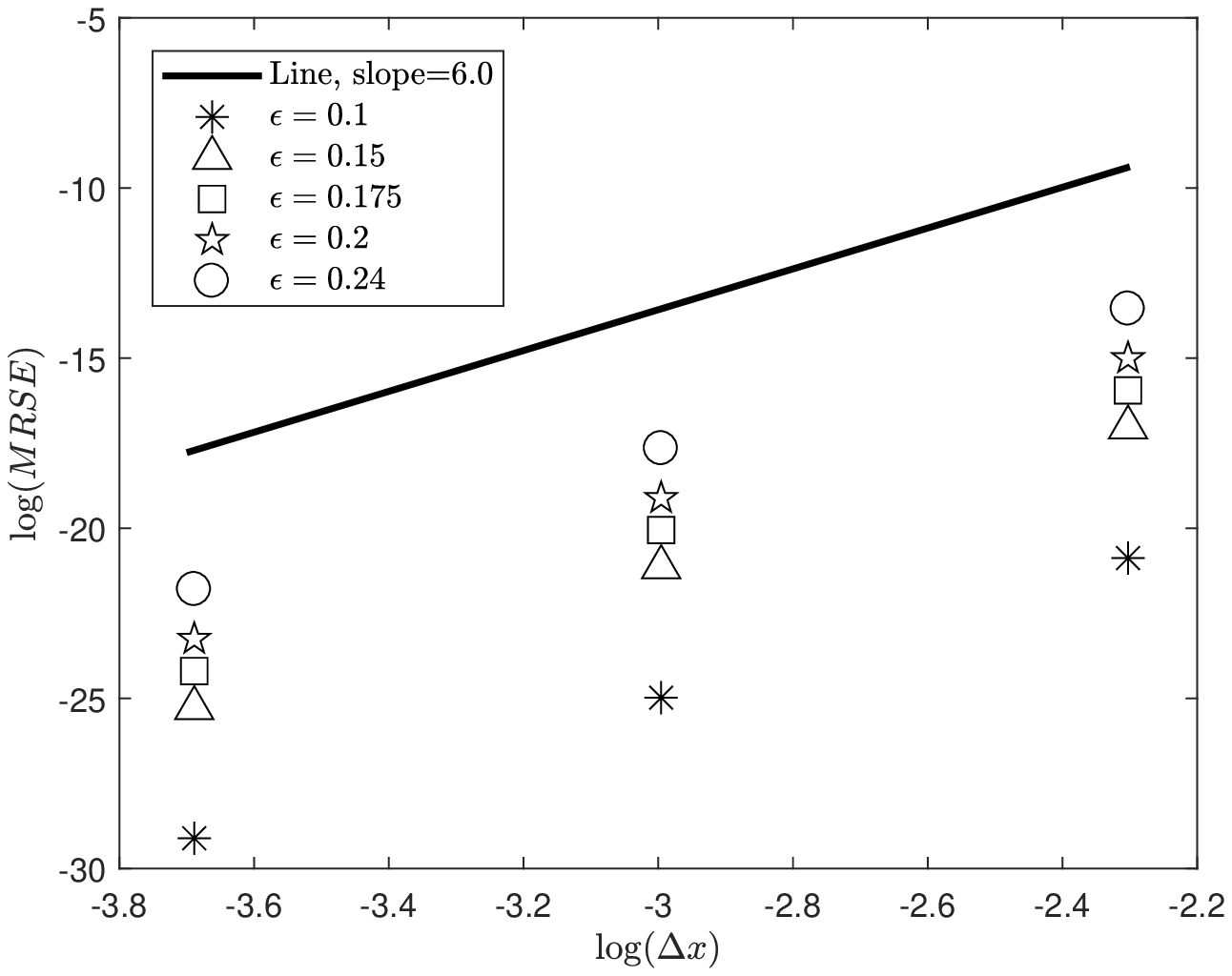}
\caption{The convergence rate of the MRT-LB model based macroscopic four-level FD scheme under the conditions of Eqs. (\ref{eq3-5}) and (\ref{eq3-6}).}
\label{fig3}
\end{figure}

\begin{figure}
\centering
\includegraphics[width=0.8\textwidth]{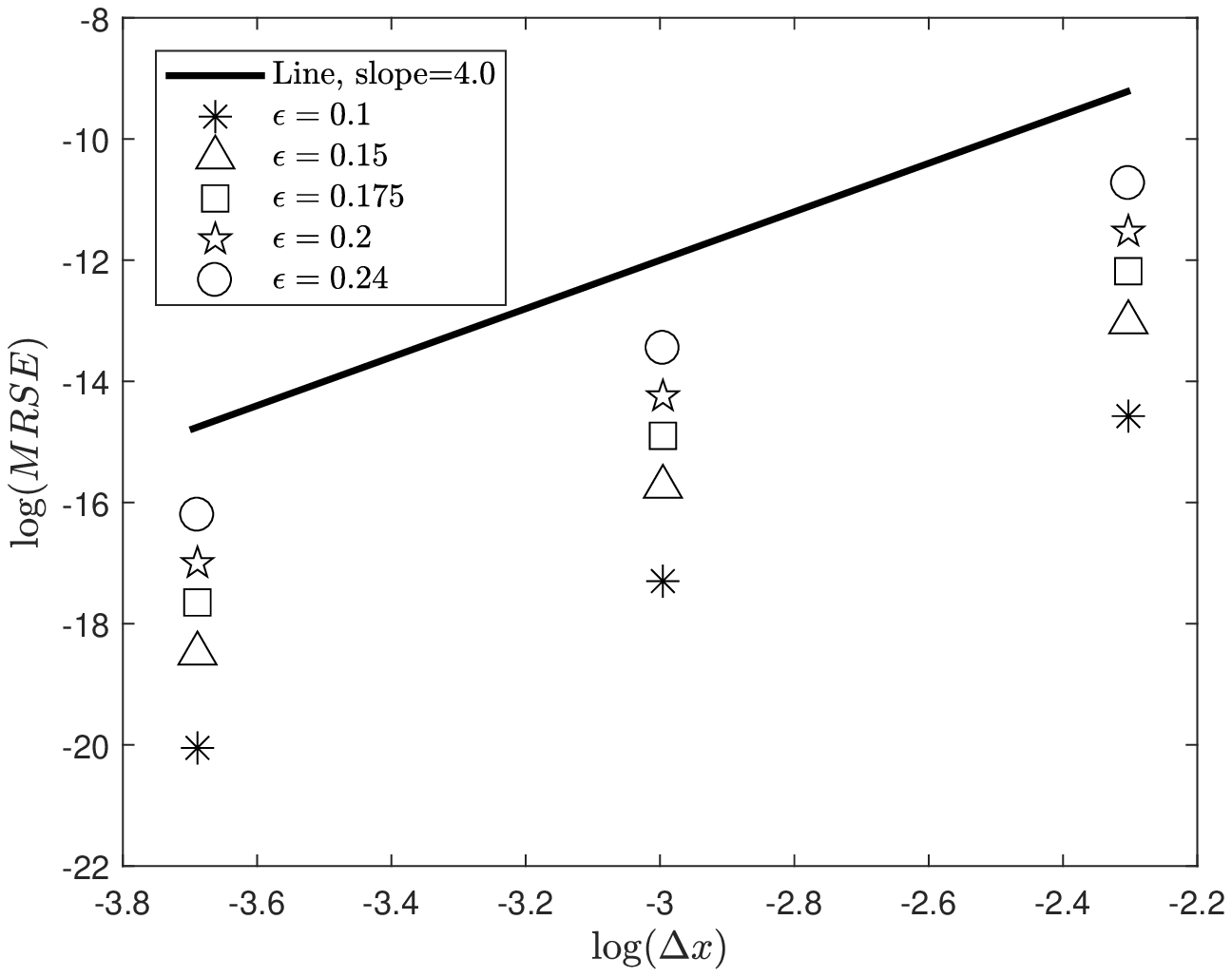}
\caption{The convergence rate of the MRT-LB model based macroscopic four-level FD scheme under the conditions of Eqs. (\ref{eq3-5}) and (\ref{eq3-6a}).}
\label{fig4}
\end{figure}

\begin{figure}
\centering
\includegraphics[width=0.8\textwidth]{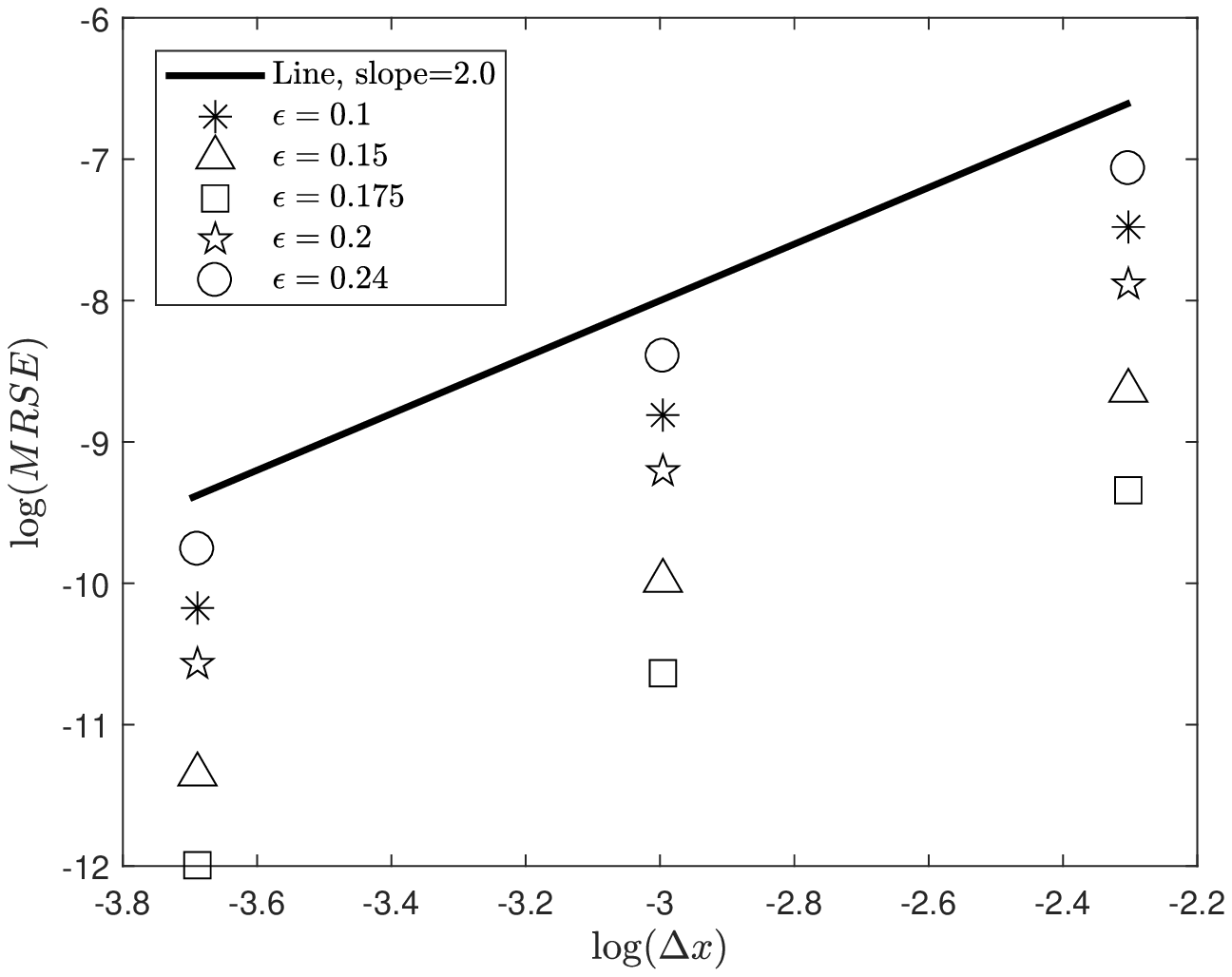}
\caption{The convergence rate of the MRT-LB model based macroscopic four-level finite-difference scheme under the condition of (\ref{eq3-5}).}
\label{fig5}
\end{figure}

We now focus on the $RMSE$ and convergence rate of the macroscopic four-level numerical scheme (\ref{eq2-21}). For this purpose, we conducted some simulations under different values of space step $\Delta x$ and mesh Fourier number $\epsilon$ [see Table 1 for details, the weight coefficient $\omega_{0}$ and relaxation parameters $s_{1}$ and $s_{2}$ are given by Eqs. (\ref{eq3-5}) and (\ref{eq3-6})], and measured the $RMSE$ between the analytical and numerical solutions at time $t=12$. As seen from Table 2 and Fig. \ref{fig3}, the MRT-LB model based macroscopic four-level numerical scheme indeed has a sixth-order convergence rate in space once the conditions of Eqs. (\ref{eq3-5}) and (\ref{eq3-6}) are satisfied, and the numerical results with a smaller mesh Fourier number $\epsilon$ are more accurate. However, as demonstrated in the previous accuracy analysis, if the conditions of Eqs. (\ref{eq3-5}) and (\ref{eq3-6}) are not met, the MRT-LB model based macroscopic numerical scheme could not achieve the sixth-order accuracy in space. To confirm this statement, we also carried out some simulations under the conditions of Eqs. (\ref{eq3-5}) and (\ref{eq3-6a}) [Eq. (\ref{eq3-6b}) is not satisfied], and presented the errors at different space steps in Table 3 and Fig. \ref{fig4}. From these table and figure, one can observe that the MRT-LB model based macroscopic four-level numerical scheme is just fourth-order accurate in space. Moreover, if both Eqs. (\ref{eq3-6a}) and (\ref{eq3-6b}) are not satisfied in our simulations, the MRT-LB model based macroscopic four-level numerical scheme, as the commonly used LB method \cite{Chai2016}, only has a second-order convergence rate, as reported in Table 4 and Fig. \ref{fig5}. We note that these results are in agreement with our theoretical analysis.
\begin{table}\label{tab_1}
\caption{The values of some parameters in the MRT-LB based macroscopic four-level numerical scheme under the conditions of Eqs. (\ref{eq3-5}) and (\ref{eq3-6}).}
\begin{center}
\begin{tabular}{c|ccc}
\hline
	$\epsilon$  & $\omega_{0}$ & $s_1$ & $s_2$  \\
  \hline
$0.1$ & $0.8310204592587027$ & $0.9159290534201945$ & $1.1450386147380731$ \\
 \hline
$0.15$ & $0.8101626131270389$ & $0.775103705680168$ & $1.1476236168426883$ \\
 \hline
$0.175$ & $0.8370678725639358$ & $0.6352970255557769$ & $1.1776696173022918$ \\
 \hline
$0.2$ & $0.870066309422671$ & $0.49037716562528605$ & $1.2047312964902426$  \\
 \hline
$0.24$ & $0.9274277013170459$ & $0.2626707812024917$ & $1.2388413217086902$ \\
\hline
\end{tabular}
\end{center}
\end{table}

\begin{table}\label{tab_2}
\caption{The $RMSE$ and $C_R$ of the MRT-LB based macroscopic four-level numerical scheme under the conditions of Eqs. (\ref{eq3-5}) and (\ref{eq3-6}).}
\begin{center}
\begin{tabular}{c|cccc}
\hline
	$\epsilon$  & $RMSE_{\Delta x=0.1}$ & $RMSE_{\Delta x=0.05}$ & $RMSE_{\Delta x=0.025}$ & $C_R$ \\
  \hline
$0.1$ & $8.59 \times 10^{-10}$ & $1.42 \times 10^{-11}$ & $2.57 \times 10^{-13}$ & $\sim5.85$ \\
 \hline
$0.15$ & $3.99 \times 10^{-8}$ & $6.56 \times 10^{-10}$ & $1.04 \times 10^{-11}$ & $\sim5.95$ \\
 \hline
$0.175$ & $1.19 \times 10^{-7}$ & $1.95 \times 10^{-9}$ & $3.11 \times 10^{-11}$ & $\sim5.95$ \\
 \hline
$0.2$ & $3.04 \times 10^{-7}$ & $5.00 \times 10^{-9}$ & $7.96 \times 10^{-11}$ & $\sim5.95$ \\
 \hline
$0.24$ & $1.31 \times 10^{-6}$ & $2.15 \times 10^{-8}$ & $3.43 \times 10^{-10}$ & $\sim5.95$ \\
\hline
\end{tabular}
\end{center}
\end{table}

\begin{table}\label{tab_3}
\caption{The $RMSE$ and $C_R$ of the MRT-LB based macroscopic four-level numerical scheme under the conditions of Eqs. (\ref{eq3-5}) and (\ref{eq3-6a}).}
\begin{center}
\begin{tabular}{c|ccccccc}
\hline
	$\epsilon$  & $\omega_{0}$ & $s_1$ & $s_2$ & $RMSE_{\Delta x=0.1}$ & $RMSE_{\Delta x=0.05}$ & $RMSE_{\Delta x=0.025}$ & $C_R$ \\
  \hline
$0.1$ & $0.8$ & $1$ & $12/11$ & $4.68 \times 10^{-7}$ & $3.08 \times 10^{-8}$ & $1.96 \times 10^{-9}$ & $\sim3.95$ \\
$0.15$ & $0.7$ & $1$ & $42/41$ & $2.21 \times 10^{-6}$ & $1.46 \times 10^{-7}$ & $9.30 \times 10^{-9}$ & $\sim3.95$ \\
$0.175$ & $0.65$ & $1$ & $78/79$ & $5.13 \times 10^{-6}$ & $3.39 \times 10^{-7}$ & $2.16 \times 10^{-8}$ & $\sim3.95$ \\
$0.2$ & $0.6$ & $1$ & $18/19$ & $9.84 \times 10^{-6}$ & $6.49 \times 10^{-7}$ & $4.14 \times 10^{-8}$ & $\sim3.95$ \\
$0.24$ & $0.52$ & $1$ & $78/89$ & $2.19 \times 10^{-5}$ & $1.44 \times 10^{-6}$ & $9.16 \times 10^{-8}$ & $\sim3.95$ \\
\hline
\end{tabular}
\end{center}
\end{table}

\begin{table}\label{tab_4}
\caption{The $RMSE$ and $C_R$ of the MRT-LB based macroscopic four-level numerical scheme under the condition of Eq. (\ref{eq3-5}) [Eqs. (\ref{eq3-6a}) and (\ref{eq3-6b}) are not satisfied].}
\begin{center}
\begin{tabular}{c|ccccccc}
\hline
	$\epsilon$  & $\omega_{0}$ & $s_1$ & $s_2$ & $RMSE_{\Delta x=0.1}$ & $RMSE_{\Delta x=0.05}$ & $RMSE_{\Delta x=0.025}$ & $C_R$ \\
  \hline
$0.1$ & $0.8$ & $1$ & $1$ & $5.65 \times 10^{-4}$ & $1.49 \times 10^{-4}$ & $3.81 \times 10^{-5}$ & $\sim1.95$ \\
$0.15$ & $0.7$ & $1$ & $1$ & $1.77 \times 10^{-4}$ & $4.62 \times 10^{-5}$ & $1.17 \times 10^{-5}$ & $\sim1.96$ \\
$0.175$ & $0.65$ & $1$ & $1$ & $8.77 \times 10^{-5}$ & $2.40 \times 10^{-5}$ & $6.18 \times 10^{-6}$ & $\sim1.91$ \\
$0.2$ & $0.6$ & $1$ & $1$ & $3.76 \times 10^{-4}$ & $1.00 \times 10^{-4}$ & $2.57 \times 10^{-5}$ & $\sim1.93$ \\
$0.24$ & $0.52$ & $1$ & $1$ & $8.55 \times 10^{-4}$ & $2.27 \times 10^{-4}$ & $5.79 \times 10^{-5}$ & $\sim1.94$ \\
\hline
\end{tabular}
\end{center}
\end{table}

\section{Conclusions}\label{section6}
In this work, we first developed a MRT-LB model for 1D diffusion equation where the D1Q3 lattice structure is considered, and then obtained a mesoscopic MRT-LB model based macroscopic four-level FD scheme. Through the theoretical analysis, one can find that the macroscopic four-level numerical scheme is unconditionally stable, and can achieve the sixth-order accuracy in space once the weight coefficient $\omega_{0}$ and the relaxation parameters $s_{1}$ and $s_{2}$ satisfy Eqs. (\ref{eq3-5}) and (\ref{eq3-6}). And also, if only the conditions of Eqs. (\ref{eq3-5}) and (\ref{eq3-6a}) are met, the macroscopic four-level numerical scheme would have a fourth-order convergence rate in space. Moreover, if only Eq. (\ref{eq3-5}) is satisfied, the macroscopic four-level numerical scheme would be second-order accurate in space, which is the same as the commonly used LB method. In addition, compared to the previous work \cite{Suga2010}, the present macroscopic numerical scheme could be more accurate through adjusting the weight coefficient $\omega_{0}$ and the relaxation parameters $s_{1}$ and $s_{2}$ properly.

We conducted some simulations to test the MRT-LB model based macroscopic four-level FD scheme, and found that the numerical results agree well with our theoretical analysis. Additionally, it is also interesting that the mesoscopic two-level LB method can be equivalent to a macroscopic four-level FD scheme, this indicates that the distribution function based mesoscopic method is more efficient than the corresponding partial-differential-equation based macroscopic numerical scheme, this may be because more information is included in the distribution function. 

Finally, we would like to point out that in the D1Q3 LB method based macroscopic numerical schemes, the present macroscopic numerical scheme (\ref{eq2-21}) with the sixth-order convergence rate is the most accurate for 1D diffusion equations. We also note that the present results are only limited to 1D case, the MRT-LB model based macroscopic numerical schemes for two and three-dimensional diffusion problems would be considered in a future work.

\bibliographystyle{siamplain}
\bibliography{references}

\end{document}